\input ./style/arxiv-general.cfg
\documentclass[seceqn,MSNbibl,number,citesort,dvips]{arxbj}
\makeatletter
   \@ifpackageloaded{graphicx}{}{\usepackage{graphicx}}
\makeatother
\aid{0}
\volume{22}
\issue{1}
\pubyear{2016}
\firstpage{325}
\lastpage{344}
\doi{10.3150/14-BEJ660} % kopijuoti is 'New paper accepted'
\docsubty{FLA}

\makeatletter

\newcommand{\afrac}[2]{#1/(#2)}

\newcommand{\rrVert}{\Vert}
\newcommand{\llVert}{\Vert}
\newtheorem{theorem}{Theorem}[section]

\newproclaim{definition}[theorem]{Definition}
\newremark{Comment}{Comment}
\newremark{Comments}{Comments}
\def\R{\mathbb{R}}
\def\N{\mathbb{N}}
\def\Z{\mathbb{Z}}
\def\E{\mathbb{E}}
\renewcommand{\P}{\mathbb{P}}
\def\Card{\operatorname{Card}}
\def\Vc{\mathcal{V}}
\def\sN{\mathcal{N}}
\def\x{\mathcal{X}}
\def\Ph{\widehat{P}}

\def\Vh{\widehat{V}}

\def\pen{\operatorname{pen}}
\newcommand{\eqref}[1]{(\ref{#1})}

\newcommand{\Ind}[1]{ \mathbf{1}_{\{#1\}}}

\makeatother

\begin{document}
\begin{frontmatter}

\title{Sharp oracle inequalities and slope heuristic for specification
probabilities estimation
in discrete random fields}
\runtitle{Oracle inequalities for specification probabilities estimation}

\begin{aug}
%%%% inicialai - be tarpu
\author[1]{\inits{M.}\fnms{Matthieu}~\snm{Lerasle}\corref{}\thanksref{1}\ead[label=e1]{mlerasle@unice.fr}} \and
\author[2]{\inits{D.Y.}\fnms{Daniel Y.}~\snm{Takahashi}\thanksref{2}\ead[label=e2]{takahashiyd@gmail.com}}
%%\runauthor{} %% auto
%\dedicated{}
\address[1]{CNRS, LJAD, UMR 7351, Univ. Nice Sophia Antipolis, 06100 Nice, France.\\ \printead{e1}}
\address[2]{Department of Psychology, Neuroscience Institute, Princeton
University, Princeton, NJ 08648, USA. \printead{e2}}
\end{aug}

% HISTORY:
\received{\smonth{12} \syear{2011}}
\revised{\smonth{6} \syear{2014}}

% ABSTRACT
%
\begin{abstract}
We study the problem of estimating the one-point specification
probabilities in non-necessary finite discrete random fields from
partially observed independent samples. Our procedures are based on
model selection by minimization of a penalized empirical criterion. The
selected estimators satisfy sharp oracle inequalities in $L_{2}$-risk.

We also obtain theoretical results on the slope heuristic for this
problem, justifying the slope algorithm to calibrate the leading
constant in the penalty. The practical performances of our methods are
investigated in two simulation studies. We illustrate the usefulness of
our approach by applying the methods to a multi-unit neuronal data from
a rat hippocampus.
\end{abstract}

% KEYWORDS
% visi is mazosios raides ir pagal abecele
%
\end{frontmatter}

%s1 #&#
\section{Introduction}
The main motivation for our work comes from neuroscience where the
advancement of multichannel and optical technology enables researchers
to record signals from tens to thousands of neurons simultaneously
\cite{Takahashi10}. The question is then to understand the interactions
between neurons in the brain and their relationships with the animal
behavior \cite{Schneidman06,Brown04}.

Following \cite{Schneidman06}, we model interactions between neurons by
discrete random fields. A discrete random field is a triplet $(S,A,P)$
where $S$ is a discrete set of \emph{sites}, possibly infinite, $A$ is
a finite alphabet, and $P$ is a probability measure on the set $\x
(S)=A^{S}$ of \emph{configurations} on $S$. Given a random field
$(S,A,P)$, we define the \emph{one point specification probabilities}
of $P$ as regular versions of the following conditional probabilities,
\[
\forall i\in S, \forall x\in\x(S),\quad\quad P_{i|S}(x)=P\bigl(x(i)|x(j), j\in S/
\{i\}\bigr).
\]
The specification probabilities are important in the applications as
they encode the conditional independence between the sites, see, for
example, \cite{Bento09,Bresler08,Csiszar06,GOT10,PW10,LT10}. The
main goal of this paper is to provide good estimators of the
specification probabilities, assuming that the configurations are only
observed on a finite subset $V_M\subset S$. Consider i.i.d. random
variables $X_{1:n}=X_{1},\ldots,X_{n}$ with common distribution $P$, the
data set is given by $(X_{i}(j))_{i=1,\ldots,n;j\in V_{M}}$. Following
\cite{BM97,BBM99,BM01}, we use a penalized criterion to select a subset
$\Vh\subset V_M$ with cardinality $\mathrm{O}(\log n)$ and show that the
empirical conditional probabilities $\Ph_{i|\Vh}$ satisfy a sharp
oracle inequality (see Section~\ref{Sect:Preliminaries} and Theorems
\ref{theo:ModelSelection} for details).

In most of the applications, the support $V_\star$ of $P_{i|S}$
(i.e., the minimal set $V_\star\subset S$ such that $P_{i|V_\star} =
P_{i|S}$) is the object of interest and the literature focus on the
estimation of $V_{\star}$, see \cite
{Bento09,Bresler08,Csiszar06,GOT10,PW10} for example. This approach
requires in general strong
assumptions on the random field, for example, it is assumed that the
data is generated by an Ising model with restrictive conditions on the
temperature parameter \cite{Bento09,GOT10,PW10}. In particular,
\cite{Bento09,Bresler08,PW10} assumed that the set $S$ is finite and that
all the sites are observed, that is, that $V_{M}=S$. When $V_{M}$ does
not contain $V_{\star}$, the meaning of the estimators in these papers
is not clear. \cite{Csiszar06} considered $S=\Z^d$ but assumed that
$V_{\star}$ is finite. Finally, \cite{GOT10,LT10} worked with infinite
sets of sites and without prior bounds on the number of interacting
sites but required a two-letters alphabet $A$ and some assumptions on
$P$ that the practitioner cannot easily verify. These restrictions are
severe in practice, for example, in neuroscience, and cast doubt on the
theoretical support for application of these methods. Our approach does
not suffer from these drawbacks. In particular, the alphabet size $|A|$
can be larger than $2$, $P$ does not need to be an Ising or Potts
model, and some configurations on $V_M$ can be forbidden. Furthermore,
$V_{\star}$ can be infinite and therefore not contained in $V_M$.

The second result of the paper is a proof of the slope heuristic for
the estimation of one-point specification probabilities in discrete
random fields. The slope heuristic was introduced in \cite{BM07} for
Gaussian model selection and has been theoretically studied only for
very few specific models \cite
{BM07,AM08,Le09,Ler2010mixing,AB10,Sa13}. Our proof technique is novel
and sheds new lights on this
phenomenon.

The paper is organized as follows. Section~\ref{Sect:Preliminaries}
presents the framework and some notations used all along the paper.
Section~\ref{Section:MainResults} introduces our estimators and the
oracle inequalities that they satisfy. In Section~\ref{sect:Bias}, the
bias for Gibbs models is computed and Section~\ref
{Section:SlopeHeuristic} is devoted to the slope heuristic.
Section~\ref{Section:Simulations} illustrates the results of previous sections
using two simulation experiments and in Section~\ref{sec:application}
our methods are applied on a neurophysiology data set.
%and Section~\ref{section:Discussion} discuss the results, making a
%detailed comparison with other works on similar problems.
The proofs of the main theorems are postponed to the Appendix~C.
%Section~\ref{Section:Proofs} and the probabilistic tools used in the
%main proofs are proved in Section~\ref{Section:Appendix}.
The methods of this article can be adapted to
%obtain oracle inequalities and the slope heuristic for
the K\"ullback loss; the interested reader can find these developments
in Section~C of the Appendix (Supplementary Materialy, \cite{supp}).
%
%s2 #&#
\section{Setting}\label{Sect:Preliminaries}
Let $(S,A,P)$ be a discrete random field, that is, a triplet where $S$
is a discrete set, $A$ is a finite set, with cardinality $|A|$ and $P$
is a probability measure on $\x(S)=A^S$. Let $V_M$ be a finite subset
of $S$ with cardinality $M\geq3$ and let $i\in S$ denote a fixed site
so that we will often omit the dependence on $i$ of some quantities
when there is no confusion. For any $x\in\x(S)$ and any $V\subset V_M$,
let $\x(V)=A^{V}$, $v=|V|$, $x(V)=(x(j))_{j\in V}$. Let $X_1,\ldots,X_n$
be i.i.d. random variables with distribution $P$. The empirical
probability measure $\Ph$ is\vspace*{1pt} defined for any $x\in\x(S)$ by
$
\Ph(x)=\frac{1}n\sum_{k=1}^n\Ind{X_k=x}$,
where $\Ind{X_k=x}=1$ if $X_k=x$ and $0$ otherwise.
The measures $P$ and $\Ph$ define probability measures on $\x(V)$ by
the formulas $P(x(V))=\int_{y\in\x(S);y(V)=x(V)}\,\mathrm{d}P(y(S))$, $\Ph
(x(V))=\sum_{y\in\x(S);y(V)=x(V)}\Ph(y)$. Hereafter, $Q$ always denotes
either $P$ or $\Ph$. For any $V\subset V_M$, $x\in\x(S)$, let
$Q_{i|V}(x)=\frac{Q(V\cup\{ i \})}{Q(V\setminus\{
i \})}$ if
$Q(V\setminus\{ i \})\ne0$, $|A|^{-1}$ otherwise. Let also
\[
P_{i|S}(x)=P\bigl(x(i) |x\bigl(S\setminus\{ i \}\bigr)\bigr)
\]
be a regular version of the conditional distribution of $P$. For any
function $f\dvtx \x(S)\to\R$, let
\[
\llVert f\rrVert_{Q}=\sqrt{\int f^2(x)
\frac{\mathrm{d}Q(x(S/\{i\}))}{|A|}}.
\]
%
%%
%For all subsets $V\subset S$, all $i \in S$, and all $x \in\x(S)$, let
%\[
%\Ph_{i|V}(x) = \frac{\Ph(x(V \cup\{i\}))}{\Ph(x(V \setminus\{i\}))},
%\]
%if $\Ph(x(V \setminus\{i\})) > 0$ and $\Ph_{i|V}(x) = |A|^{-1}$,
%otherwise.
The observation set is $X_{1:n}(V_M)=(X_1(j),\ldots,X_n(j))_{j\in V_M}$.
Algebraic computations show
\[
\forall y\in\x(V_M),\quad\quad \Ph(y)=\frac{1}n\sum
_{i=1}^n\Ind{X_i(V_M)=y},
\]
and for any $V\subset V_M$, $\Ph(x(V))=\sum_{y\in\x
(V_M);y(V)=x(V)}\Ph
(x(V))$ can be computed from the data set. Hence, for $V \subset V_M$
the empirical probability $\Ph_{i|V}$ is an estimator of $P_{i|S}$. The
\emph{$L_{2,P}$-risk} of $\Ph_{i|V}$ is defined by $\llVert\Ph
_{i|V}-P_{i|S}\rrVert^2_{P}$. We can decompose the risk via Pythogoras
relation (see Proposition B.11)
\[
\llVert\Ph_{i|V}-P_{i|S}\rrVert_{P}^2=
\llVert\Ph_{i|V}-P_{i|V}\rrVert_{P}^2+
\llVert P_{i|V}-P_{i|S}\rrVert_{P}^2.
\]
The random term $\llVert\Ph_{i|V}-P_{i|V}\rrVert_{P}^2$ is called the
\emph{variance} and the deterministic term $\llVert
P_{i|V}-P_{i|S}\rrVert_{P}^2$ is called the \emph{bias}. Let $s\geq
3$ be an integer and let
\[
\Vc_s= \{ V\subset V_M,v\leq s \},\quad\quad
N_s=\Card( \Vc_s ).
\]
An \emph{oracle} is a set $V_{o} \in\Vc_s$ that minimizes the risk,
\textit{that is},
\[
\llVert\Ph_{i|V_o}-P_{i|S}\rrVert_{P}^2
= \min_{V \in\Vc_s}\llVert\Ph_{i|V}-P_{i|S}\rrVert
_{P}^2
\]
and the minimal risk is called \emph{oracle risk}. We will show in the
next section that we can obtain an estimator $\Vh$ such that the risk
of $\Ph_{i|\Vh}$ is close to the oracle risk.
%
%s3 #&#
\section{Model selection results}\label{Section:MainResults}
Let start with a concentration inequality for the variance term of the risks.

%
%th3.1 #&#
\begin{theorem}\label{theo:controlVar}
Let $Q\in\{ P,\Ph \}$ and let $V\in\Vc_s$. Then, for
all $\delta>1$
and all $0<\eta\leq1$,
%
%e3.1 #&#
\begin{equation}
\label{eq:controlVarP} \P\biggl(\llVert\Ph_{i|V}-P_{i|V}\rrVert
_{Q}^2>\frac{6}{|A|} \biggl((1+8\eta)
\frac{|A|^{v}}n+\frac{4\log(2\delta)}{\eta n}+\frac{9\log
(2\delta
)^2}{\eta^{4} n} \biggr) \biggr)\leq
\frac{1}{\delta}.
\end{equation}
%
%\begin{equation}\label{eq:controlVarPn}
%\left\|\Ph_{i|V}-P_{i|V}\right\|^{2}_{\Ph}\leq\frac{6}{|A|}\left((1+8
%\eta)\frac{|A|^{v}}n+\frac{4\log(2\delta)}{\eta n}+\frac{9\log(2
%\delta)^2}{\eta^{4} n}\right).
%\end{equation}
\end{theorem}

\begin{Comment*}
%\begin{itemize}
%\item
The bound can be integrated to give the following control
\[
\E\bigl[\llVert\Ph_{i|V}-P_{i|S}\rrVert_{P}^2\bigr]=
\llVert P_{i|V}-P_{i|S}\rrVert_{P}^2+C
\frac{|A|^{v-1}}{n},
\]
for some absolute constant $C$. This control depends on the
approximation properties of $V$ through the bias $\llVert
P_{i|V}-P_{i|S}\rrVert_{P}^2$ and on the variance via the upper bound
$|A|^{v-1}/n$. Our goal now is to find a subset $V$ that balances these
two terms. This is precisely the aim of the following result.
\end{Comment*}
%\end{itemize}

%
%th3.2 #&#
\begin{theorem}\label{theo:ModelSelection}
Let
\[
\Vh=\arg\min_{V\in\Vc_s} \bigl\{-\llVert\Ph_{i|V}\rrVert
_{\Ph
}^2+\pen(V) \bigr\},\qquad\mbox{where }\pen(V)\geq12
\frac{|A|^{v-1}}{n}.
\]
There exists a constant $\kappa=\kappa(|A|)$ such that, with
probability larger than $1-\delta^{-1}$,
%
%e3.2 #&#
\begin{equation}
\label{eq:BetterOrIneqq} \llVert P_{i|S}-\Ph_{i|\Vh}\rrVert
_{P}^2\leq\biggl(1+\frac
{8}{\log
(\delta)} \biggr)\inf
_{V\in\Vc_s} \bigl\{\llVert P_{i|S}-P_{i|V}\rrVert
_{P}^2+\pen(V) \bigr\}+\kappa\frac{(\log(N_s^2\delta))^2}{n}.
\end{equation}
\end{theorem}

\begin{Comments*}
\begin{itemize}
\item The bound can be integrated and yields
\[
\E\bigl[\llVert P_{i|S}-\Ph_{i|\Vh}\rrVert_{P}^2\bigr]
\leq C_1\inf_{V\in
\Vc
_s} \biggl\{\llVert
P_{i|S}-P_{i|V}\rrVert_{P}^2+
\frac
{|A|^{v-1}}{n} \biggr\}+C_2\frac{(s\log M)^2}n,
\]
for some absolute constant $C_1$ and a constant $C_2$ depending only on
$|A|$. Therefore, $\Vh$ optimizes the bound given by Theorem~\ref
{theo:controlVar}, up to the residual $(s\log(M))^2$ term, among all
the subsets of $\Vc_s$.
\item Enlarging the number of observed sites makes the control over all
subsets in $\Vc_s$ harder, leading to a $(s\log M)^2$ loss in the
rates. On the other hand, it is helpful to reduce the bias as will be
shown in the next section.
\item A very interesting feature of this result for the applications is
that it holds without restrictions on $P$ and the size of $A$ or $S$ in
$(S,A,P)$.
\end{itemize}
\end{Comments*}

%
%s4 #&#
\section{Computation of the bias}\label{sect:Bias}
To complete the study of our estimator, it remains to understand the
bias
$\llVert P_{i|S}-P_{i|V}\rrVert_{P}^2$. We present two important
examples where explicit upper bounds can be obtained.
%s4.1 #&#
\subsection{The Ising model}
Let $S=\Z^d$ and let $(J_{i,j})_{(i,j)\in S^2}$ be an interaction
potential, which is a collection of real numbers such that for any
$i\ne j\in S$, $J_{i,i}=0$, $J_{i,j}=J_{j,i}$ and
\[
\beta:= \sup_{i\in S}\sum_{j\in S}
| J_{i,j} |<\infty.
\]
The parameter $1/\beta$ is also called the temperature parameter in the
physic literature where the model was initially introduced, see \cite
{Giorgii88}. The Ising model is the triplet $(S,A,P)$, where $A=
\{ -1,1 \}$ and $P$ is given by its specifications by
\[
P_{i|S}(x)=\frac{\mathrm{e}^{\sum_{j\in S}J_{i,j}x(i)x(j)}}{\mathrm
{e}^{\sum_{j\in
S}J_{i,j}x(i)x(j)}+\mathrm{e}^{-\sum_{j\in S}J_{i,j}x(i)x(j)}}=\frac
{1}{1+\mathrm{e}^{-2\sum
_{j\in S}J_{i,j}x(i)x(j)}} .
\]
It follows from Theorem~4.5 in \cite{LT10} that
\[
\llVert P_{i|S}-P_{i|V} \rrVert_{P}\leq\sup
_{x\in\x(S)}\bigl|P_{i|S}(x)-P_{i|V}(x)\bigr|\leq
C_\beta\sum_{j\notin V}| J_{i,j} |
.
\]
Rates of convergence can be obtained from this bound and our model
selection theorem. For example, let $d_{\infty}(i,j) = \max\{
|i_k-j_k|\dvt k \in\{1, \ldots, d\}\}$, assume that $s\log M=\mathrm
{O}((\log
n)^2)$ and that there exists constants $r$ and $r'$ such that $\sum
_{j\in S: d_\infty(i,j)>k}| J_{i,j} |\leq k^{-r}$ and $\sum
_{j>k}| J_{i,j}^{*} |\leq\mathrm{e}^{-r'k}$, where
$J^*_{i,j}$ denote the
rearrangement of the $J_{i,j}$ by decreasing absolute values. Then, for
any $i\in V_M$, denoting by $\alpha_i$ the largest real number such
that $\{j \in\Z\dvt d_\infty(i,j) \leq n^{\alpha_i}\} \subset V_M$, we have
\begin{eqnarray*}
\E\bigl[\llVert P_{i|S}-P_{i|\hat{V}} \rrVert_P^2\bigr]
&\leq& C\frac{(\log
n)^4}n+C_\beta\bigl( n^{-\alpha_i r}+
n^{-\afrac{2r'}{2r'+\log2}} \bigr)
\\
&\leq& C_{\beta} n^{- ( \alpha_i r\wedge\afrac{2r'}{2r'+\log2}
)} .
\end{eqnarray*}
Other consequences of this bound obtained under different assumptions
on the $(J_{i,j})_{i,j\in S}$ are discussed in Section~\ref{sec:ProofForBias}.

%s4.2 #&#
\subsection{The Gibbs model}
Assume that $A$ is a finite set of real numbers in $[-1,1]$, $S=\Z^d$
for some $d\geq1$. Let $ ( (J^{(k)}_{i,i_1,\ldots
,i_k})_{(i,i_1,\ldots,i_k)\in S^{k+1}} )_{k\geq0}\in\prod_{k\geq
0}\R
^{k+1}$ be a collection of real numbers such that
\[
\sum_{k\geq0}\sum_{(i,i_1,\ldots,i_k)\in S^{k+1}}
\bigl| J^{(k)}_{i,i_1,\ldots,i_k} \bigr|=\beta<\infty.
\]
For any $x\in\x(S)$ and $i\in S$, denote by
\[
J_i(x)=\sum_{k\geq0}\sum
_{(i_1,\ldots,i_k)\in
S^{k}}J^{(k)}_{i,i_1,\ldots,i_k}\prod
_{\ell=1}^kx(i_\ell) .
\]
Suppose that the conditional probabilities can be written in the
following way:
\[
P_{i|S}(x)=\frac{\mathrm{e}^{x(i)J_i(x)}}{\sum_{a\in A}\mathrm
{e}^{aJ_i(x)}} .
\]
The triplet $(S,A,P)$ is called a \emph{Gibbs} model,
Ising models are special instances of Gibbs models
where for all $k\ge2$ and all $(j_1,\ldots,j_k)\in S^k$,
$J_{i,j_1,\ldots,j_k}=0$. For any $\ell\leq M$, denote by
$(J^*_{i,\ell,n})_{n=1,\ldots,M^\ell}$ the rearrangement of
the $J^{(\ell)}_{i,i_1,\ldots,i_\ell}$ by decreasing absolute values.
We consider the following assumption.\def\theequation{J}
%
%e4.1 #&#
\begin{equation}
\label{eq:Cond3J} \forall\ell,n\in\N^*,\quad\quad  \sum_{r\ge n}
\bigl| J^*_{i,\ell ,r} \bigr|\leq\beta\mathrm{e}^{-\gamma\ell
^{2+\alpha} n} ,
\end{equation}
for some constant $\gamma$ and $\alpha>0$.
Under Assumption \eqref{eq:Cond3J}, we can build a set $V$ with
cardinality $v\leq\frac{1+2\alpha}{\gamma\alpha+\log|A|(1+2\alpha
)}\log n$ such that
the bias term is upper bounded \def\theequation{\arabic{section}.\arabic{equation}}
by\setcounter{equation}{0}
%
%e4.1 #&#
\begin{equation}
\label{eq:ContBiasGibbs} \llVert P_{i|S}-P_{i|V} \rrVert
_P^2\leq C_{\alpha,\beta,\gamma,|A|} \biggl( \frac{(\log
n)^{1/(2+\alpha)}}{n^{\afrac{\alpha\gamma}
{\gamma\alpha+\log|A|(1+2\alpha)} }}+
\sum_{\ell\ge1}\sum_{i_1,
\ldots i_\ell\in S:\exists j; i_j\notin V_M}
\bigl| J^{(\ell)}_{i,i_1, \ldots,i_\ell} \bigr| \biggr) .
\end{equation}
The bound \eqref{eq:ContBiasGibbs} is proved in Section~\ref
{sec:ProofForBias}. From Theorem~\ref{theo:controlVar} and $v\leq
\frac
{1+2\alpha}{\gamma\alpha+\log|A|(1+2\alpha)}\log n$, for some absolute
constant $C$,
\begin{eqnarray*}
\E\bigl[ \llVert\Ph_{i|V}-P_{i|V} \rrVert_P^2\bigr]
\leq C\frac
{|A|^{v-1}}n=\frac
{C}{|A|n^{\afrac{\alpha\gamma}{\gamma\alpha+\log|A|(1+2\alpha)}
}} .
\end{eqnarray*}
Therefore, for some constant $C_{\alpha,\beta,\gamma,|A|}$ and rate
$\theta=\frac{2\alpha\gamma}{2\alpha\gamma+(1+2\alpha)\log|A|}$,
\begin{eqnarray*}
\E\bigl[ \llVert P_{i|S}-\Ph_{i|\Vh} \rrVert_P^2\bigr]
\leq C_{\alpha
,\beta,\gamma
,|A|} \biggl[ \biggl( \frac{\log n}{n} \biggr)^{\theta}+
\sum_{\ell
\leq v}\sum_{i_1,\ldots i_\ell\in S:\exists j; i_j\notin O}
\bigl| J^{(\ell)}_{i,i_1,\ldots,i_\ell} \bigr| \biggr] .
\end{eqnarray*}

%s5 #&#
\section{Slope heuristic}\label{Section:SlopeHeuristic}
%In practice, the constant $6K|A|^{-1}$ in Theorem
%\ref{theo:ModelSelection} is a bit pessimistic. In order to optimize
%this constant,
The slope heuristic was introduced in \cite{BM07}. Let
%
%e5.1 #&#
\begin{equation}
\label{eq:defVh} \Vh=\arg\min_{V\in\Vc_{s}} \bigl\{-\llVert
\Ph_{i|V}\rrVert_{\Ph}^2+\pen(V) \bigr\}.
\end{equation}
The heuristic states that there exist a minimal penalty $\pen_{\min}$
and a complexity measure (to be defined) satisfying the following properties.
\begin{enumerate}[SH3]
\item[SH1] When $\pen(V)<(1-\eta)\pen_{\min}(V)$, the complexity
of $\Vh
$ is as large as possible.
\item[SH2] When $\pen(V)=(1+\eta)\pen_{\min}(V)$, the complexity
of $\Vh
$ is much smaller.
\item[SH3] When $\pen(V)=2\pen_{\min}(V)$, the risk of $\Vh$ is
equivalent to the oracle risk.
\end{enumerate}
%
%The heuristic is used to calibrate the constant in front of the
%penalty and the choice of measure of complexity depends on the problem.
The purpose of this section is to justify this heuristic. We will show
some theoretical evidence for the slope heuristic using $\Delta
_V=\llVert \Ph_{i|V}-P_{i|V} \rrVert^2_{\Ph}$ as a complexity
measure for $V$ and as a
minimal penalty. It may be useful for the intuition to make the
following approximation $n\Delta_V/|A|^v\approx C$ although it is only
proved in Theorem~\ref{theo:controlVar} that $\E[\Delta_V]\le
C|A|^{v}/n$. For example, this explains why it's natural to consider
$\Delta_V$ as a measure of complexity.
The following theorem gives some theoretical grounds justifying SH1.

%
%th5.1 #&#
\begin{theorem}\label{theo:MinimalPenalty}
Let $r>0$, $\epsilon>0$. Let $\Vh$ be defined by \eqref{eq:defVh} and
assume that
\[
\P\bigl(\forall V\in\Vc_{s},0\leq\pen(V)\leq(1-r)\llVert\Ph
_{i|V}-P_{i|V}\rrVert_{\Ph}^2 \bigr)
\geq1-\epsilon.
\]
Then, for all $\delta>2$, with probability larger than $1-\epsilon
-2\delta^{-1}$,
\[
\llVert P_{i|\Vh}-\Ph_{i|\Vh}\rrVert_{\Ph}^2
\geq\sup_{V\in\Vc
_{s}} \bigl\{r\llVert P_{i|V}-
\Ph_{i|V}\rrVert_{\Ph}^2-2\llVert
P_{i|S}-P_{i|V}\rrVert_{P}^2 \bigr\}-
\frac{17}3\frac{(\log
(N_s^2\delta))^2}n.
\]
\end{theorem}

\begin{Comments*}
\begin{itemize}
\item Let us give some intuition on this result. Algebraic
computations, see \eqref{eq:MinimCrit}, show that $\Vh$ minimizes, up
to centered remainder terms, the quantity
%
%e5.2 #&#
\begin{equation}
\label{eq:approx.crit} \llVert P_{i|S}-P_{i|V}\rrVert
_{P}^2+\pen(V)-\llVert\Ph_{i|V}-P_{i|V}
\rrVert_{\Ph}^2 .
\end{equation}
We assume in Theorem~\ref{theo:MinimalPenalty} that $\pen(V)=(1-\eta
)\Delta_V$, thus $\Vh$ minimizes the bias minus $\eta\Delta_V$. When
the bias term decreases with $V$, as in the models presented in
Section~\ref{sect:Bias} and when $n\Delta_V/|A|^v\approx C$, both terms
decrease with $V$ and the minimum is achieved for $\Vh=V_M$. Thus,
$\Vh
$ maximizes the complexity $\Delta_V$.
\item Theorem~\ref{theo:MinimalPenalty} makes this statement more
precise, showing that this result actually holds when, for $V=V_M$,
both the bias and the logarithmic remainder term are negligible
compared to the variance part of the risk.
%Theorem~\ref{theo:MinimalPenalty} is therefore a minimal penalty
%theorem. It states that, if the penalty term is too small, the
%complexity of the selected model (measured here with $\left\|P_{i|V}-
%\Ph_{i|V}\right\|_{\Ph}^2$) is as large as possible.
\end{itemize}
\end{Comments*}

Let us now turn to the associated optimal penalty theorem which proves
SH2 and SH3.

%th5.2 #&#
\begin{theorem}\label{theo:OptimalPenalty}
Let $\delta>5 $, $r_2\ge r_1>0$, $\epsilon>0$ and assume that
%
%e5.3 #&#
\begin{eqnarray}
\label{HypOptPen} \P\biggl(\forall V\in\Vc_{s},(1+r_1)
\leq\frac{\pen(V)}{\llVert
\Ph
_{i|V}-P_{i|V}\rrVert_{\Ph}^2}\leq(1+r_2) \biggr)\geq1-\epsilon.
\end{eqnarray}
Let $\Vh$ be defined by \eqref{eq:defVh}. For all $V$ in $\Vc_{s}$, let
$p_-^{V}=\inf_{x\in\x(V),P(x(V))\neq0}P(x(V))$ and assume that, for
some $\varepsilon\leq1$,
\[
\inf_{V\in\Vc_s}p_-^{V}\geq\varepsilon^{-2}
\frac{\log
(nN_s\delta)}{n}.
\]
Then, there exists an absolute constant $C$ such that, with probability
larger than $1-5\delta^{-1}-\epsilon$, for all $V$ in $\Vc_{s}$, for
all $\eta>0$,
%
%e5.4 #&#
\begin{eqnarray}
\label{ineqOracleOptimal} \frac{(1-\eta)\wedge
(r_1-C(1+r_1)\varepsilon)}{(1+\eta)\vee
(r_2+C(1+r_2)\varepsilon)} \llVert P_{i|S}-\Ph_{i|\Vh}
\rrVert_{P}^2\leq\llVert P_{i|S}-
\Ph_{i|V}\rrVert_{P}^2+\frac{6}{\eta}
\frac{(\log
(N_s^2\delta))^2}{n}.
\end{eqnarray}
\end{theorem}

\begin{Comments*}
\begin{itemize}
\item In this theorem, following \cite{AM08}, the main task is to show that
%
%e5.5 #&#
\begin{equation}
\label{eq:todoforslope} \Delta_V\simeq\llVert\Ph_{i|V}-P_{i|V}
\rrVert_{P}^2.
\end{equation}
When \eqref{HypOptPen} holds with $r_1=r_2=r$, then
\[
\pen(V)=(1+r)\Delta_V\simeq\Delta_V+r\llVert\Ph
_{i|V}-P_{i|V}\rrVert_{P}^2.
\]
From \eqref{eq:approx.crit}, $\Vh$ minimizes the sum of the bias and
$r$ times the variance. The complexity should thus be much smaller,
which proves SH2 for $\pen_{\min}(V)=\Delta_V$. Theorem~\ref
{theo:OptimalPenalty} shows that the complexity of the selected model,
that is bounded by the risk, is actually upper bounded by the supremum
between the oracle risk and the remainder term, at least when
$\varepsilon$ is small enough.

%\item Choose first $r_1=r_2>0$, $\eta=1/2$, the penalty is equal to
%$(1+r_1)\pen_{\min}(V)$. It comes from (\ref{ineqOracleOptimal}) that,
%if $\varepsilon\rightarrow0$,
%\[
%C_{r_1}\Ex{\left\|P_{i|\Vh}-\Ph_{i|\Vh}\right\|_{P}^2}\leq\inf_{V\in
%\Vc_s}\Ex{\left\|P_{i|S}-\Ph_{i|V}\right\|_{P}^2}+\frac{(s\log
%M)^2}{n}.
%\]
%Moreover, from \eqref{eq:todoforslope} (see Lemma~\ref{lem:Ph-P}), $
%\Ex{\left\|P_{i|\Vh}-\Ph_{i|\Vh}\right\|_{\Ph}^2}$ satisfies the same
%inequality. Under the hypotheses of the previous remark, the
%complexity of $\Vh$ is upper bounded by the oracle risk $\inf_{V
%\subset V_M}\Ex{\norma{P_{i|S}-\Ph_{i|V}}^2}$ and should be much
%smaller than the biggest one. Hence, we observe a jump of the
%complexity of the selected model around $\pen_{\min}$, this is SH2.
%
\item Take then $r_1=r_2=1$, that is, a penalty equal to
\[
\pen(V)=2\pen_{\min}(V)\simeq\Delta_V+\llVert\Ph
_{i|V}-P_{i|V}\rrVert_{P}^2.
\]
Then \eqref{eq:approx.crit} shows that $\Vh$ minimizes an approximately
optimal criterion, and $\Ph_{i|\Vh}$ satisfies an oracle inequality
that is asymptotically optimal, which proves SH3. Inequality (\ref
{ineqOracleOptimal}) makes this result more precise, showing that the
oracle inequality is indeed asymptotically optimal when the oracle rate
of convergence is larger than the remainder term. Moreover, in this
case, the rate of convergence of the leading quantity in the oracle is
driven by the supremum of the rates $\eta$ and $\varepsilon$.
% We have therefore justified the slope heuristic for the $L_2$-risk.
%In the following section, we give the theorems justifying it for the K
%\"ullback loss.
\end{itemize}
\end{Comments*}

Theorem~\ref{theo:OptimalPenalty} cannot be used directly to build an
estimator since the complexity is unknown. Nevertheless, Theorem~\ref
{theo:controlVar} shows that $\Delta_V$ is upper bounded by $K\Theta
_V$, with $\Theta_V=|A|^{v-1}/n$ and some constant $K$ that may not be
optimal. This suggests to consider penalties of the form $K\Theta_V$,
for some $K$ that has to be optimized. To achieve this goal, \cite
{AM08} proposed the following algorithm.
\begin{enumerate}
\item For all $K>0$, denote by $\Vh(K)$ the model selected with $\pen
(V)=K\Theta_V$.
\item Find $K_{\min}$ such that $\Theta_{\Vh(K)}$ is very large for
$K<K_{\min}$ and much smaller for $K>K_{\min}$.
\item Select $\Vh=\Vh(2K_{\min})$.
\end{enumerate}
This algorithm is based on the slope heuristic. Indeed, assume that
$\pen_{\min}(V)=K_0\Theta_V$ for some unknown $K_0$. Then, $K_{\min}$
shall be close to $K_0$ because we observe a jump of the complexity
$\Theta_{\Vh}$ around $K_{\min}\Theta_V$ as expected by SH1, SH2.
Therefore, $\Vh$, chosen by $2K_{\min}\Theta_V\simeq2\pen_{\min}(V)$
shall be optimal from SH3. We did not prove that this algorithm
improves the choice of $K$ in theory but the simulation study of the
next section presents examples where it does in practice.

%s6 #&#
\section{Simulation studies} \label{Section:Simulations}

In this section, we illustrate the results obtained in previous ones
using simulation experiments. All the simulations were implemented by a
set of MATLAB\textsuperscript{\textregistered} routines that can be
downloaded from \href
{http://www.princeton.edu/\textasciitilde
dtakahas/publications/LT11routines.zip}{www.princeton.edu/\textasciitilde
dtakahas/publications/LT11routines.zip}.

\begin{figure}

\includegraphics{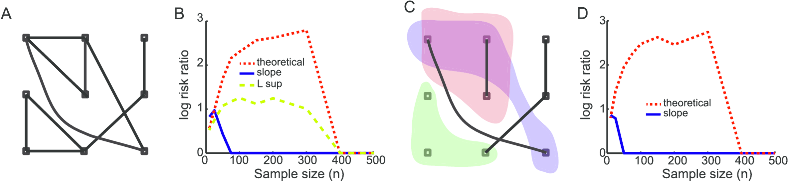}

\caption{Simulation study. (A)~Representation of the interacting pairs
of the Ising model used in the first simulation experiment.
The numbering of the sites increases from the top left to the
bottom right. (B)~Performance of the model selection for the first
experiment. Plot of the log risk ratio for the model selection
procedure
using $K=2$ (dotted red line), optimizing the constant using the slope
heuristic (solid blue line), using the $L_{\infty}$-risk method with
slope
heuristic (dashed yellow). (C)~Representation of the interacting
neurons
of the Gibbs model used in the second simulation experiment. The
colored
regions represent the three-way interactions. (D)~Performance of the
model
selection for the second experiment. The legend is the same as in~(B).}\label{fig:NewSimulationRiskRatio}
\end{figure}

Let $S=\{1, \ldots, 9\}$ and $A = \{-1,1\}$. For the first simulation,
we consider an Ising model $(S,A,P)$, with one-point specification
probabilities given by
\[
\forall x \in\x(S), \quad\quad P_{i|S}(x) = \frac{1}{1+\exp(-2\sum_{j
\in
S}J_{ij}x(i)x(j))},
\]
where the $J_{ij}$'s are given by $J_{1,2} = J_{1,5} = -J_{2,5} =
J_{1,9} = J_{2,9} = J_{3,6} = -J_{4,7} = -J_{4,8} = -J_{7,8} = J_{6,8}
= 0.5$. The rest of $J_{ij}$'s are equal to zero. For each $i \in S$,
the pair of sites $(i,j)$ where $j \in V_i$ is shown in Figure~\ref
{fig:NewSimulationRiskRatio}(A). For the first experiment, we study the
site $i = 9$ and its interaction sites. We simulate independent samples
of the Ising model and compare the performances of the model selection
procedures given by (1) the penalty given in Theorem~\ref
{theo:ModelSelection} (theoretical), (2) the same penalty, but using
the slope algorithm described in Section~\ref{Section:SlopeHeuristic}
to calibrate the constant in front of $|A|^{v-1}/n$, and (3) the
$L_\infty$-risk method with slope heuristic proposed in~\cite{LT10}.
The performances of the estimators are measured by the logarithm of the
ratio between the risk of the estimated model and the oracle risk.
Figure~\ref{fig:NewSimulationRiskRatio}(B) shows the median value of the
risk ratio calculated for $100$ independent replicas. The maximum
number of allowed interacting sites was set to $s=5$. The simulations
were done for increasing sample sizes $n=10$, $25$, $50$, $75$, $100$,
$150$, $200$, $300$, $400$, $500$.

For the second simulation, we consider a Gibbs model $(S,A,P)$, with
one-point conditional probabilities given by
\[
\forall x \in\x(S),\quad\quad P_{i|S}(x) = \frac{1}{1+\exp(-2\sum_{j
\in
S}J_{ij}x(i)x(j) + \sum_{k \in S}\sum_{j \in S}J_{ijk}x(i)x(j)x(k))}.
\]
The non-null pairwise interactions are given by $ -J_{2,5} = J_{1,9} =
J_{3,6} = J_{6,8} = 0.5$, and the three-way interactions are specified
by $J_{1,2,5} = J_{1,2,9} = -J_{4,7,8} = 0.5$. The rest of $J_{ij}$'s
and $J_{ijk}$'s are equal to zero. For each $i$, the interacting
neighborhood $V_i$ is shown in Figure~\ref
{fig:NewSimulationRiskRatio}(C). We show the results for $i = 9$. We
compute the risk ratio as in the first experiment (Figure~\ref
{fig:NewSimulationRiskRatio}(D)). The simulations are done for increasing
sample sizes $n=10$, $25$, $50$, $75$, $100$, $150$, $200$, $300$,
$400$, $500$ (Figure~\ref{fig:NewSimulationRiskRatio}(D)). Observe that
in both experiments the slope heuristic improves the performance of the
model selection, allowing to recover the oracle even for data set as
small as $50$ in our examples. For this example, any method that uses
the Ising model to estimate the parameters has a non-null bias and
therefore the risk will be strictly larger than the oracle risk.
Further simulations are shown in the Appendix (Section~C).

%f1 #&#

%s7 #&#
\section{Application to multi-unit neuronal data} \label{sec:application}

In this section, we illustrate the usefulness of the proposed methods
on experimental data set. In neuroscience, it is conjectured that the
set of interacting neurons represents different animal behaviors \cite
{Schneidman06}. Modifications of the graph of interacting neurons for
different tasks have been repeatedly shown \cite{Schneidman06}.
Nevertheless, if this hypothesis has any validity, we expect the set of
interacting neurons to be the same when the same task is performed. We
used our method here to test this hypothesis, which seems to be less
verified in the literature.

The data set used contains multichannel simultaneous recordings made
from layer CA1 of the right dorsal hippocampus of a Long-Evans rat
during open field tasks in which the animal chased randomly placed
drops of water while on a elevated square platform. It was downloaded
from \href
{http://crcns.org/data-sets/hc/hc-2/about-hc-2}{http://crcns.org/data-sets/hc/hc-2/about-hc-2}.
Details about the recording technique and experimental set up can be
found at the website or in \cite{Mizuseki09}.

The spiking data set used is ec016.430.res.1, ec016.430.res.2,
ec016.430.res.3, ec016.430.  res.4, ec016.430.res.5, ec016.430.res.6,
ec016.430.res.7, ec016.430.res.8. The full data set contains a total of
55 isolated neurons. For the analysis, we kept only the $11$ neurons
that showed more than $30\,000$ spikes during the experiment. The data
set was sampled at 20~kHz. We binned the data with non-overlapping bins
of size $10$~ms. If there was at least one spike in the bin, we coded it
as $+1$, otherwise we coded as $-1$.
The spiking activity of the 11 neurons was recorded for 106.8 minutes.
To ensure independence of the observations, we subsampled the data
using one observation at each $500$~ms, which is an order of magnitude
larger than a typical decay of correlation (when the correlation
becomes zero) between neurons in time. We then splitted the data into
two parts, one sample for the first half of the experiment ($n=64\,099$,
first $53.4$~min) and another sample for the second half of the
experiment ($n=64\,099$, second $53.4$~min).

\begin{figure}[t]

\includegraphics{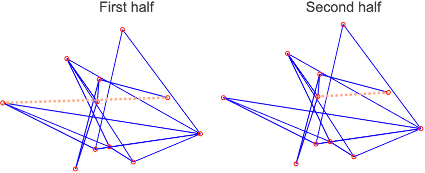}

\caption{Representation of the interacting neuronal sites for
the first half and second half of the experiment. The edges
between sites indicate the interacting pairs. The dotted orange
edges indicate the interactions that differed between both conditions.
Observe that the interactions are represented by a graph for
convenience
of visualization, but for our method the interactions are not
restricted
to pairwise interaction as shown by our theoretical results and in
Figure~\protect\ref{fig:NewSimulationRiskRatio}(D).}\label
{fig:ExperimentNeuronSite}
\end{figure}

We computed our estimators of the interacting neurons and calibrate the
constant in front of the penalty with the slope algorithm described in
the end of Section~\ref{Section:SlopeHeuristic}. For each site, the
maximum number of allowed interacting sites was $s=3$.
Figure~\ref{fig:ExperimentNeuronSite} shows the results obtained for
the first and second parts of the experiment. We clearly see that the
interacting neuronal sites remained stable, with only one pair of
interaction that changed between the two data sets. This result,
together with those in the literature showing changes in interacting
neighborhoods for different behaviors, corroborates the hypothesis that
the set of interacting neurons can be related to specific animal behavior.

%f2 #&#

\begin{appendix}
%%%%%%%%%%%%%%%%%%%%%%%%%%%%%%%%%%%%%%%%%%%%%%%%%%%%%%%%%%%%
%sA #&#
\section{Proofs}\label{Section:Proofs}
%sA.1 #&#
\subsection{Proof of Theorem \texorpdfstring{\protect\ref{theo:controlVar}}{3.1}}
Let $\theta>0$ to be chosen later and let $Q$ denote either $P$ or
$\Ph
$. We decompose the risk as follows
\begin{eqnarray*}
\llVert\Ph_{i|V}-P_{i|V}\rrVert_{Q}^2&=&
\sum_{x\in\x(V)}\frac
{Q(x(V/\{
i\}))}{|A|} \bigl(
\Ph_{i|V}(x)-P_{i|V}(x) \bigr)^2
\\
&=&\sum_{x\in\x(V),Q(x(V/\{i\}))\leq\theta(|A|^{v}n)^{-1}}\frac
{Q(x(V/\{i\}))}{|A|} \bigl(
\Ph_{i|V}(x)-P_{i|V}(x) \bigr)^2
\\
&&{}+\sum_{x\in\x(V),Q(x(V/\{i\}))> \theta(|A|^{v}n)^{-1}}\frac
{Q(x(V/\{
i\}))}{|A|} \bigl(
\Ph_{i|V}(x)-P_{i|V}(x) \bigr)^2.
\end{eqnarray*}
As the cardinal of $\x(V)$ is $|A|^{v}$ and $ (\Ph
_{i|V}(x)-P_{i|V}(x) )^2\leq1$, the first term in this
decomposition is upper bounded by $\theta n^{-1}$. Hence
%
%eA.1 #&#
\begin{equation}
\label{eq:interm1} \llVert\Ph_{i|V}-P_{i|V}\rrVert
_{Q}^2=\frac{\theta}n+\sum
_{x\in
\x(V),
Q(x(V/\{i\}))> \theta(|A|^{v}n)^{-1}}\frac{Q(x(V/\{i\}))}{|A|} (\Ph
_{i|V}-P_{i|V}
)^2.
\end{equation}
Hereafter in the proof of Theorem~\ref{theo:controlVar}, we denote by
\[
\x^{\theta}(V)= \bigl\{x\in\x(V)\dvt  Q\bigl(x\bigl(V/\{i\}\bigr
)\bigr)>
\theta\bigl(|A|^{v}n\bigr)^{-1} \bigr\}.
\]
It comes from Lemma B.1 that
\begin{eqnarray*}
&&\llVert\Ph_{i|V} -P_{i|V}\rrVert_{P}^2-
\frac{\theta
}n
\\
&&\quad=\sum_{x\in\x^{\theta}(V)}\frac{P(x(V/\{i\}))}{|A|}
\bigl(\Ph_{i|V}(x)-P_{i|V}(x) \bigr)^2
\\
&& \quad\leq \sum_{x\in\x^{\theta}(V)}\frac{ (|\Ph
(x(V))-P(x(V))|+\Ph_{i|V}(x)| (\Ph(x(V/\{i\}))-P(x(V/\{i\}))
)| )^2}{|A|P(x(V/\{i\}))}
\\
&& \quad\leq \frac{2}{|A|} \biggl( \sum_{x\in\x^{\theta}(V)}
\frac{
(\Ph(x(V))-P(x(V)) )^2}{P(x(V/\{i\}))}+\sum_{x\in\x^{\theta
}(V/\{ i\} )}\frac{ (\Ph(x(V/\{i\}))-P(x(V/\{i\}))
)^2}{P(x(V/\{i\}))}
\biggr).
\end{eqnarray*}
From Lemma B.1, we also have
\[
\bigl|\Ph_{i|V}(x)-P_{i|V}(x)\bigr|\leq\frac{|\Ph(x(V))-P(x(V))|
+P_{i|V}(x)| (\Ph(x(V/\{i\}))-P(x(V/\{i\})) )|}{|A|\Ph
(x(V/\{i\}))}.
\]
Hence
\begin{eqnarray*}
&&\bigl|\Ph_{i|V}(x)-P_{i|V}(x)\bigr|
\\
&&\quad\leq\frac{|\Ph(x(V))-P(x(V))|
+(P_{i|V}(x)+\Ph_{i|V}(x))| (\Ph(x(V/\{i\}))-P(x(V/\{i\}
)) )|}{|A|\sqrt{\Ph(x(V/\{i\}))P(x(V/\{i\}))}}.
\end{eqnarray*}
Thus,
\[
\llVert\Ph_{i|V} -P_{i|V}\rrVert_{\Ph}^2-
\frac
{\theta
}n=\sum_{x\in\x^{\theta}(V)}\frac{\Ph(x(V/\{i\}))}{|A|}
\bigl(\Ph_{i|V}(x)-P_{i|V}(x) \bigr)^2
\]
is smaller than
\begin{eqnarray*}
&&\sum_{x\in\x^{\theta}(V)}\frac{ (|\Ph
(x(V))-P(x(V))|+(\Ph_{i|V}(x)+P_{i|V}(x))| (\Ph(x(V/\{i\}
))-P(x(V/\{i\}
)) )| )^2}{|A|P(x(V/\{i\}))}
\\
&&\quad\leq\frac{2}{|A|} \biggl( \sum_{x\in\x^{\theta}(V)}
\frac{
(\Ph(x(V))-P(x(V)) )^2}{P(x(V/\{i\}))}+2\sum_{x\in\x^{\theta
}(V/\{ i\} )}\frac{ (\Ph(x(V/\{i\}))-P(x(V/\{i\}))
)^2}{P(x(V/\{i\}))}
\biggr).
\end{eqnarray*}
We use Theorem B.8 with $b=\sqrt{\theta^{-1} |A|^v n}$,
for all $x>0$, for all $\eta>0$, we have, with probability larger than
$1-2\mathrm{e}^{-x}$,
\[
\llVert\Ph_{i|V}-P_{i|V}\rrVert_{Q}^2
\leq\frac{\theta}n+\frac
{6}{|A|} \biggl((1+\eta)^3
\frac{|A|^{v}}n+\frac{4x}{\eta n}+\frac{32|A|^v
x^2}{\theta\eta^3 n} \biggr).
\]
Take $\theta=8|A|^{v/2}x\eta^{-3/2}$, we obtain
\[
\llVert\Ph_{i|V}-P_{i|V}\rrVert_{Q}^2
\leq\frac{6}{|A|} \biggl((1+\eta)^3\frac{|A|^{v}}n+
\frac{4x}{\eta n}+\frac{6|A|^{v/2} x}{\eta^{3/2}
n} \biggr).
\]
Using $ab\leq\eta a^2+(4\eta)^{-1}b^2$, we finally get
\[
\llVert\Ph_{i|V}-P_{i|V}\rrVert_{Q}^2
\leq\frac{6}{|A|} \biggl((1+8\eta)\frac{|A|^{v}}n+\frac{4x}{\eta n}+
\frac{9x^2}{\eta^{4} n} \biggr).
\]
%

%

%
%sA.2 #&#
\subsection{Proof of Theorem \texorpdfstring{\protect\ref{theo:ModelSelection}}{3.2}}
The theorem follows from the slightly more general following result.

%
%thA.1 #&#
\begin{theorem}
Let $K>1$ and let
\[
\Vh=\arg\min_{V\in\Vc_s} \bigl\{-\llVert\Ph_{i|V}\rrVert
_{\Ph
}^2+\pen(V) \bigr\},\quad\quad\mbox{where }\pen(V)\geq6K
\frac{|A|^{v-1}}{n}.
\]
Then, there exists a constant $\kappa=\kappa(|A|,K)$ such that for all
$\delta\ge1$, with probability larger than $1-\delta^{-1}$,
%
%eA.2 #&#
\begin{equation}
\label{eq:ClassOrIneq} \llVert P_{i|S}-\Ph_{i|\Vh}\rrVert
_{P}^2\leq\kappa\biggl(\inf_{V\in\Vc
_s}
\bigl\{\llVert P_{i|S}-P_{i|V}\rrVert_{P}^2+
\pen(V) \bigr\} +\frac
{(\log(N_s^2\delta))^2}{n} \biggr).
\end{equation}
Moreover, when $K\ge2$, there exists a constant $\kappa=\kappa(|A|,K)$
such that, with probability larger than $1-\delta^{-1}$,
%
%eA.3 #&#
\begin{equation}
\label{eq:BetterOrIneq} \llVert P_{i|S}-\Ph_{i|\Vh}\rrVert
_{P}^2\leq\biggl(1+\frac
{8}{\log
(\delta)} \biggr)\inf
_{V\in\Vc_s} \bigl\{\llVert P_{i|S}-P_{i|V}\rrVert
_{P}^2+\pen(V) \bigr\}+\kappa\frac{(\log(N_s^2\delta))^2}{n}.
\end{equation}
\end{theorem}

\begin{pf}
For $Q\in\{ P,\Ph \}$, let $(\cdot,\cdot)_{Q}$ be the scalar
product associated
to the $L_{2,Q}$-norm $\|\cdot\|_{Q}$. Let $V$ and $V'$ in the collection
$\Vc_{s}$. We have
\begin{eqnarray*}
&&\frac{1}{|A|}\sum_{x\in\x(V\cup V')} \Ph\bigl(x\bigl(V\cup
V'\bigr)\bigr)P_{i|V}(x)
\\
 &&\quad=\sum_{x\in\x(V)}\frac{\Ph(x(V/\{i\}))}{|A|} \Ph
_{i|V}(x)P_{i|V}(x)= (\Ph_{i|V},P_{i|V}
)_{\Ph},
\\
&&\frac{1}{|A|}\sum_{x\in\x(V\cup V')} P\bigl(x\bigl(V\cup
V'\bigr)\bigr)P_{i|V}(x)=\sum
_{x\in\x
(V)}\frac{P(x(V/\{i\}))}{|A|} P^2_{i|V}(x)=
\llVert P_{i|V}\rrVert^2_{P}.
\end{eqnarray*}
Hence, for all $V$, $V'$ in $\Vc_{s}$,
%
%eA.4 #&#
\begin{eqnarray}\label{eq:DecompositionCrit}
\llVert\Ph_{i|V}\rrVert_{\Ph}^2&=&\llVert
P_{i|V}\rrVert_{\Ph
}^2+2 (\Ph_{i|V}-P_{i|V},P_{i|V}
)_{\Ph}+\llVert\Ph_{i|V}-P_{i|V}\rrVert
_{\Ph}^2
\nonumber
\\
&=&\llVert P_{i|V}\rrVert_{P}^2+\llVert
\Ph_{i|V}-P_{i|V}\rrVert_{\Ph
}^2- \bigl(
\llVert P_{i|V}\rrVert_{\Ph}^2-\llVert
P_{i|V}\rrVert_{P}^2 \bigr)
\\
&&{}+\frac{2}{|A|}\sum_{x\in\x(V\cup V')}\bigl(\Ph\bigl(x
\bigl(V\cup V'\bigr)\bigr)-P\bigl(x\bigl(V\cup V'
\bigr)\bigr)\bigr)P_{i|V}(x).
\nonumber
\end{eqnarray}
Moreover, from Pythagoras relation see Proposition B.11, we have
\[
\llVert P_{i|S}-P_{i|V}\rrVert_{P}^2=
\llVert P_{i|S}\rrVert_{P}^2-\llVert
P_{i|V}\rrVert_{P}^2.
\]
By definition of $\Vh$, we have, for all $V$ in $\Vc_{s}$,
\[
\llVert P_{i|S}\rrVert_{P}^2-\llVert
\Ph_{i|\Vh}\rrVert_{\Ph
}^2+\pen(\Vh)\leq\llVert
P_{i|S}\rrVert_{P}^2-\llVert\Ph_{i|V}
\rrVert_{\Ph
}^2+\pen(V).
\]
Hence, for all $0<\nu\leq1$, from (\ref{eq:DecompositionCrit}),
\[
\nu\llVert P_{i|S}-\Ph_{i|\Vh}\rrVert_{P}^2
\leq\llVert P_{i|S}-P_{i|\Vh}\rrVert_{P}^2+
\nu\llVert P_{i|\Vh}-\Ph_{i|\Vh
}\rrVert_{P}^2
\]
is smaller than
%
%eA.5 #&#
\begin{eqnarray}\label{eq:DecompositionRisk}
&&\llVert P_{i|S}-P_{i|V}\rrVert_{P}^2+
\pen(V)-\llVert\Ph_{i|V}-P_{i|V}\rrVert_{\Ph}^2\nonumber
\\
&&\quad{}-
\bigl(\pen(\Vh)-\llVert\Ph_{i|\Vh
}-P_{i|\Vh}\rrVert
_{\Ph}^2-\nu\llVert\Ph_{i|\Vh}-P_{i|\Vh
}
\rrVert_{P}^2 \bigr)
\nonumber
\\[-8pt]\\[-8pt]
&&\quad{}+ \bigl(\llVert P_{i|V}\rrVert_{\Ph}^2-\llVert
P_{i|V}\rrVert_{P}^2-\llVert P_{i|\Vh}
\rrVert_{\Ph}^2+\llVert P_{i|\Vh}\rrVert
_{P}^2 \bigr)\nonumber
\\
&&\quad{}+\frac{2}{|A|}\sum_{x\in\x(V\cup\Vh)}\bigl(\Ph\bigl(x(V
\cup\Vh)\bigr)-P\bigl(x(V\cup\Vh)\bigr)\bigr) \bigl(P_{i|\Vh
}(x)-P_{i|V}(x)
\bigr).\nonumber
\end{eqnarray}
We have also,
\begin{eqnarray*}
&&\llVert P_{i|V}\rrVert_{\Ph}^2- \llVert
P_{i|V}\rrVert_{P}^2-\llVert P_{i|\Vh}
\rrVert_{\Ph}^2+\llVert P_{i|\Vh}\rrVert
_{P}^2
\\
&&\quad=\frac{1}{|A|}\sum_{x\in\x((V\cup\Vh))}\bigl(\Ph\bigl(x
\bigl((V\cup\Vh)/\{i\} \bigr)\bigr)-P\bigl(x\bigl((V\cup\Vh)/\{
i\}\bigr)\bigr)
\bigr) \bigl(P^2_{i|V}(x)-P^2_{i|\Vh}(x)
\bigr).
\end{eqnarray*}
Let $0<\eta\leq1$, $\delta>1$ and assume that, $N_s\geq2$. Let
$\Omega
^{\delta}$ be the intersection of the following events:
%
%eA.6 #&#
%eA.7 #&#
\begin{eqnarray}
 \Omega_1^{\delta}&=& \biggl\{\forall V\in\Vc_{s},
\llVert\Ph_{i|V}-P_{i|V}\rrVert^{2}_{\Ph}
\leq\frac{6}{|A|} \biggl((1+8\eta)\frac
{|A|^{v}}n+\frac{13\log(2N_s\delta)^2}{\eta^{4} n}
\biggr) \biggr\},\nonumber
\\
\Omega_2^{\delta}&=& \biggl\{\forall V\in\Vc_{s},
\llVert\Ph_{i|V}-P_{i|V}\rrVert^{2}_{P}
\leq\frac{6}{|A|} \biggl((1+8\eta)\frac
{|A|^{v}}n+\frac{13\log(2N_s\delta)^2}{\eta^{4} n}
\biggr) \biggr\},\nonumber
\\
\Omega_3^{\delta}&=&\biggl\{\forall V,V' \in
\Vc_{s}^2, \llVert P_{i|V}\rrVert
_{\Ph}^2-\llVert P_{i|V}\rrVert
_{P}^2 -\llVert P_{i|V'}\rrVert
_{\Ph}^2+\llVert P_{i|V'}\rrVert
_{P}^2
\nonumber
\\[-8pt]\label{Def:Omega3}\\[-8pt]
&&\hphantom{\biggl\{\forall V,V' \in
\Vc_{s}^2,\quad}\leq2\llVert P_{i|V}-P_{i|V'}\rrVert_{P}\sqrt
{2\frac{\log
(N_s^2\delta)}{n}}+\frac{\log(N_s^2\delta)}{3n}\biggr\},\nonumber
\\
\Omega_4^{\delta}&=&\biggl\{\forall V,V' \in
\Vc_{s}^2, \sum_{x\in
\x
(V\cup V')}\bigl(\Ph
\bigl(x\bigl(V\cup V'\bigr)\bigr)-P\bigl(x\bigl(V\cup
V'\bigr)\bigr)\bigr)\frac
{P_{i|V'}(x)-P_{i|V}(x)}{a}
\nonumber
\\[-8pt]\label{Def:Omega4}\\[-8pt]
&&\hphantom{\biggl\{\forall V,V' \in
\Vc_{s}^2, \sum_{x\in
\x
(V\cup V')}\quad} \leq\llVert P_{i|V}-P_{i|V'}\rrVert_{P}\sqrt
{2\frac{\log
(N^2_s\delta)}{n}}+\frac{\log(N_s^2\delta)}{3n}\biggr\}.\nonumber
\end{eqnarray}
Theorem~\ref{theo:controlVar}, Lemma B.10 and union bounds give that
\[
P \bigl( \bigl(\Omega^{\delta} \bigr)^c \bigr) \leq
\frac
{4}{\delta}.
\]
For all $V$, $V'$ in $\Vc_{s}$ and all $\xi>0$, on $\Omega^{\delta}$,
we have
\begin{eqnarray*}
&&2\sum_{x\in\x(V\cup V')}\bigl(\Ph\bigl(x\bigl(V\cup V'
\bigr)\bigr)-P\bigl(x\bigl(V\cup V'\bigr)\bigr)\bigr)\frac
{P_{i|V'}(x)-P_{i|V}(x)}{|A|}+
\llVert P_{i|V}\rrVert_{\Ph}^2
\\
&&\quad{}-\llVert P_{i|V}\rrVert_{P}^2-\llVert
P_{i|V'}\rrVert_{\Ph
}^2+\llVert P_{i|V'}
\rrVert_{P}^2\leq\frac{\xi}2\llVert
P_{i|V}-P_{i|V'}\rrVert_{P}^2+ \biggl(
\frac{16}{\xi}+1 \biggr)\frac{\log(N_s^2\delta)}{3n}.
\end{eqnarray*}
From (\ref{eq:DecompositionRisk}), we deduce that, on $\Omega^{\delta
}$, for all $0<\xi<\eta$,
\begin{eqnarray*}
(\nu-\xi)\llVert P_{i|S}-\Ph_{i|\Vh}\rrVert_{P}^2&
\leq&(1+\xi)\llVert P_{i|S}-P_{i|V}\rrVert_{P}^2+
\pen(V)
\\
&&{}- \biggl(\pen(\Vh)-(1+\nu) (1+\eta)^3\frac{6}{|A|}
\frac
{|A|^{\widehat{v}}}n \biggr)
\\
&&{}+\frac{1}n \biggl(\frac{78(1+\nu)}{\eta^4|A|}\bigl(\log
(2N_s\delta
)\bigr)^2+ \biggl(\frac{16}{\xi}+1 \biggr)\log
\bigl(N_s^2\delta\bigr) \biggr).
\end{eqnarray*}
Take at first $0<\xi<\nu$ and $0<\eta$ sufficiently small to ensure
that $(1+\nu)(1+\eta)^3\leq K$ to obtain (\ref{eq:ClassOrIneq}). To
obtain (\ref{eq:BetterOrIneq}), choose $\nu=1$ and $\eta>0$
sufficiently small to ensure that $(1+\eta)^3<K/2$ and $\xi=(\log
(N_s^{2}\delta))^{-1}$. We conclude the proof, saying that the
inequality is obvious when $\delta<4$, and, when $\delta\geq4$,
\[
\frac{1+(\log N_s^2\delta)^{-1}}{1-(\log N_s^2\delta)^{-1}}=1+\frac
{2(\log N_s^2\delta)^{-1}}{1-(\log N_s^2\delta)^{-1}}\leq1+\frac
{2(\log\delta)^{-1}}{1-(\log\delta)^{-1}}\leq1+
\frac{8}{\log\delta}.
\]
\upqed\end{pf}

%

%
%%
%sA.3 #&#
\subsection{Proof of the bias control}\label{sec:ProofForBias}
%sA.3.1 #&#
\subsubsection{Discussion on the Ising model}
In this section, we discuss some consequences of the bound given on the
bias term in the Ising model, under additional assumptions on the $J_{i,j}'s$.
\begin{enumerate}
\item Assume that the set of $j\in S$ such that $J_{i,j}\neq0$, $\sN
_i$ is finite and that $\sN_i\subset V_M$. The bound \eqref
{eq:BetterOrIneq} implies that, when $\log_2(n)\geq|\sN_i|$,
\[
\E\bigl[\llVert P_{i|S}-\Ph_{i|\hat{V}} \rrVert_P^2\bigr]
\leq C\frac
{(\log(n)\log
(M))^2}n
\\
+C_\beta\frac{2^{|\sN_i|}}{n}\leq C_{\beta,|\sN_i|}\frac{(\log
(n)\log
(M))^2}n .
\]
\item Assume that there exist constants $r$ and $r'$ such that
$M=n^{r}$ and, for any $k\in\N$, $\sum_{j>k}| J^*_{i,j} |
\leq
\mathrm{e}^{-r'k}$, then
\begin{eqnarray*}
\E\bigl[ \llVert P_{i|S}-\Ph_{i|\hat{V}} \rrVert_P^2\bigr]
&\leq& Cr^2\frac
{\log
(n)^4}n+C_\beta\biggl( \biggl( \sum
_{j\notin V_M}| J_{i,j} | \biggr)^2+
n^{-\afrac{2r'}{2r'+\log2}} \biggr)
\\
&\leq& C_{r,\beta} \biggl( \biggl( \sum_{j\notin V_M}
| J_{i,j} | \biggr)^2+ n^{-\afrac{2r'}{2r'+\log2}} \biggr) .
\end{eqnarray*}
\end{enumerate}
%
%sA.3.2 #&#
\subsubsection{Proof of the bound on the bias in the Gibbs case}
%Notice that this assumption naturally extend the one made in the Ising
%model.
%Actually,
%%\eqref{eq:Cond1J} holds with $r=\infty$, \eqref{eq:Cond2J} reduces to
%$\sum_{j\in S: d_\infty(i,j)>k}\absj{J_{i,j}}\leq k^{-r'}$ in the
%Ising model and
%\eqref{eq:Cond3J} can be rewritten $\sum_{j>k}\absj{J_{i,j}^{*}}\leq
%\beta\mathrm{e}^{-\gamma k}$ in the Ising case.
%
In order to bound the bias term $\llVert P_{i|S}-P_{i|V} \rrVert
_P^2$, we still
use the inequalities
\[
\llVert P_{i|S}-P_{i|V} \rrVert_P\leq\llVert
P_{i|S}-P_{i|V} \rrVert_\infty
\leq\sup_{x,y\in\x(S):x(V\cup\{ i \})=y(V\cup
\{ i \})}\bigl| P_{i|S}(x)-P_{i|S}(y) \bigr|
.
\]
Now, we will build an approximation set $V=\bigcup_{\ell=0}^{\log
_{|A|}n}\sN_{\ell}$ and bound the bias of $P_{i|V}$, using the
inequality for any $v\leq|V|$,
\begin{eqnarray*}
&&\frac{| J_i(x)-J_i(y) |}2
\\
&&\quad\leq\sum_{\ell\leq v}\sum
_{i_1,\ldots
i_\ell\in S:\exists j; i_j\notin\sN_\ell}\bigl| J^{(\ell
)}_{i,i_1,\ldots ,i_\ell} \bigr|+\sup
_{z\in\x(S)}\sum_{\ell>v}
\bigl| J_i^{(\ell)}(z) \bigr|
\\
&&\quad\leq\sum_{\ell\leq v}\sum_{i_1,\ldots i_\ell\in S:\exists j;
i_j\notin
V_M}
\bigl| J^{(\ell)}_{i,i_1,\ldots,i_\ell} \bigr|+\sum_{i_1,\ldots
i_\ell\in
V_M:\exists j; i_j\notin\sN_\ell}
\bigl| J^{(\ell)}_{i,i_1,\ldots ,i_\ell } \bigr|+\frac{\beta
}{1-\mathrm{e}^{-\gamma}} \mathrm{e}^{-rv^{2+\alpha}} .
\end{eqnarray*}
Let $\sN_\ell$ denote the union of the $K_\ell$ $\ell$-tuples
$i_1,\ldots,i_\ell$ such that $(J^*_{i,\ell,r})_{r=1,\ldots,K_\ell}$
are indexed by the $ \{ (i,i_1,\ldots,i_\ell),\mbox{s.t.}
(i_1,\ldots,i_\ell)\in\sN_\ell \}$. $\sN_\ell$ has a
cardinality smaller than
$K_\ell\ell$ and by assumption~\eqref{eq:Cond3J}, we have
%
%eA.8 #&#
\begin{equation}
\label{eq:TBias1} \sum_{i_1,\ldots i_\ell\in O:\exists j; i_j\notin
V_\ell}\bigl| J^{(\ell )}_{i,i_1,\ldots,i_\ell} \bigr|
\leq\beta\mathrm{e}^{-\gamma\ell^{2+\alpha
}K_\ell
} .
\end{equation}
Now, let us fix some $\nu>0$ and let $K_\ell=1+ \lfloor\nu
\ell^{-2-\alpha}\log n \rfloor$ for any $\ell\leq(\nu\log
n)^{1/(2+\alpha)}$ and $K_\ell=0$
when $\ell>(\nu\log n)^{1/(2+\alpha)}$. In particular, $K_\ell\ge
\nu
\ell^{-2-\alpha}\log n$ when $\ell\leq(\nu\log n)^{1/(2+\alpha)}$,
hence, from \eqref{eq:TBias1}, for any $1\le\ell\leq(\nu\log
n)^{1/(2+\alpha)}$, we have
\[
\sum_{i_1,\ldots i_\ell\in V_M:\exists j; i_j\notin\sN_\ell}\bigl|
J^{(\ell)}_{i,i_1,\ldots,i_\ell} \bigr|
\leq\frac{\beta}{(1-\mathrm{e}^{-\gamma
})n^{\nu\gamma}} .
\]
Therefore, the bias term is upper bounded by
\[
\llVert P_{i|S}-P_{i|V} \rrVert_P^2
\leq C_{\alpha,\beta,\gamma
,|A|} \biggl( \frac{\log n}{n^{\nu\gamma}}+\sum
_{\ell\ge1}\sum_{i_1,\ldots i_\ell\in S:\exists j; i_j\notin O}
\bigl| J^{(\ell )}_{i,i_1,\ldots,i_\ell} \bigr| \biggr) .
\]
Moreover, $V$ has cardinality upper bounded by
%$|V|\leq\frac{1+2\alpha}{\gamma\alpha+\log|A|(1+2\alpha)}\log n$ such
%that
%
\[
\sum_{\ell= 1}^{(\nu\log n)^{1/(2+\alpha)}}\ell K_\ell\leq
\sum_{\ell=
1}^{(\nu\log n)^{1/(2+\alpha)}} \biggl( \ell+
\frac{\nu\log n}{\ell
^{1+\alpha}} \biggr)\leq\frac{1+2\alpha}{\alpha}\nu\log n .
\]

%
%sA.4 #&#
\subsection{Proof of Theorem \texorpdfstring{\protect\ref{theo:MinimalPenalty}}{5.1}}
Let us introduce, for all $V$ in $\Vc_{s}$,
\[
L(V)=\llVert P_{i|V}\rrVert_{P}^2-\llVert
P_{i|V}\rrVert_{\Ph
}^2+\frac
{2}{|A|}\sum
_{x\in\x(V)}\bigl(\Ph\bigl(x(V)\bigr)-P\bigl(x(V)\bigr)
\bigr)P_{i|V}(x).
\]
By definition of $\Vh$, we have, for all $V$ in $\Vc_{s}$,
\[
\llVert P_{i|S}\rrVert_{P}^2-\llVert
\Ph_{i|\Vh}\rrVert_{\Ph
}^2+\pen(\Vh)\leq\llVert
P_{i|S}\rrVert_{P}^2-\llVert\Ph_{i|V}
\rrVert_{\Ph
}^2+\pen(V).
\]
Hence from inequality (\ref{eq:DecompositionCrit}) in the proof of
Theorem~\ref{theo:ModelSelection}, we have, for all $V$ in $\Vc_{s}$,
%
%eA.9 #&#
\begin{eqnarray}\label{eq:MinimCrit}
&&\llVert P_{i|S}-P_{i|\Vh}\rrVert_{P}^2 +
\bigl(\pen(\Vh)-\llVert\Ph_{i|\Vh}-P_{i|\Vh}\rrVert
_{\Ph}^2 \bigr)-L(\Vh)
\nonumber
\\[-8pt]\\[-8pt]
&&\quad\leq\llVert P_{i|S}-P_{i|V}\rrVert_{P}^2+
\bigl(\pen(V)-\llVert\Ph_{i|V}-P_{i|V}\rrVert
_{\Ph}^2 \bigr)-L(V).\nonumber
\end{eqnarray}
Let $\Omega_{\pen}= \{0\leq\pen(V)\leq(1-r)\llVert\Ph
_{i|V}-P_{i|V}\rrVert_{\Ph}^2 \}$ and let $\Omega^{\delta
}_{\min\,
\pen}=\Omega_3^{\delta}\cap\Omega_4^{\delta}\cap\Omega_{\pen
}$, where
$\Omega_3^{\delta}$ and $\Omega_4^{\delta}$ are respectively
defined in
(\ref{Def:Omega3}) and (\ref{Def:Omega4}). It comes from Lemma B.10 and our assumption on $\pen(V)$ that
$P((\Omega^{\delta}_{\min\,\pen})^c)\leq\epsilon+2\delta^{-1}$.
Moreover, on $\Omega^{\delta}_{\min\,\pen}$, we have, for all $\eta>0$,
\begin{eqnarray*}
&&\bigl|L(\Vh)-L(V)\bigr|
\\
&&\quad\leq\eta\llVert P_{i|S}-P_{i|\Vh}\rrVert
_{P}^2+\eta\llVert P_{i|S}-P_{i|V}
\rrVert_{P}^2+ \biggl(\frac{16}{\eta}+1 \biggr)
\frac
{\log
(N_s^2\delta)}{3n},
\\
&&(1-\eta)\llVert P_{i|S}-P_{i|\Vh}\rrVert_{P}^2 -
\llVert\Ph_{i|\Vh
}-P_{i|\Vh}\rrVert_{\Ph}^2
\\
&&\quad\leq(1+\eta)\llVert P_{i|S}-P_{i|V}\rrVert
_{P}^2-r\llVert\Ph_{i|V}-P_{i|V}
\rrVert_{\Ph}^2+ \biggl(\frac{16}{\eta}+1 \biggr)
\frac{\log
(N_s^2\delta)}{3n}.
\end{eqnarray*}
We conclude the proof choosing $\eta=1$.

%

%
%sA.5 #&#
\subsection{Proof of Theorem \texorpdfstring{\protect\ref{theo:OptimalPenalty}}{5.2}}
Let
\[
\Omega_{\pen}= \bigl\{\forall V\in\Vc_{s},(1+r_1)
\llVert\Ph_{i|V}-P_{i|V}\rrVert_{\Ph}^2
\leq\pen(V)\leq(1+r_2)\llVert\Ph_{i|V}-P_{i|V}
\rrVert_{\Ph}^2 \bigr\},
\]
let $\Omega^{\delta}_{\rm comp}=\Omega_3^{\delta}\cap\Omega
_4^{\delta
}\cap\Omega_{\pen}$, where $\Omega_3^{\delta}$ and $\Omega
_4^{\delta}$
are respectively defined in (\ref{Def:Omega3}) and (\ref{Def:Omega4}).
It comes from Lemma B.10 and our
assumption on $\pen(V)$ that $P((\Omega^{\delta}_{\min\,\pen
})^c)\leq
\epsilon+2\delta^{-1}$. Moreover, on $\Omega^{\delta}_{\min\,\pen
}$, we
have, from (\ref{eq:MinimCrit}), for all $\eta>0$,
\begin{eqnarray*}
&&(1-\eta) \llVert P_{i|S}-P_{i|\Vh}\rrVert_{P}^2+r_1
\llVert\Ph_{i|\Vh
}-P_{i|\Vh}\rrVert_{P}^2+(1+r_1)
\bigl( \llVert\Ph_{i|\Vh
}-P_{i|\Vh}\rrVert_{\Ph}^2-
\llVert\Ph_{i|\Vh}-P_{i|\Vh}\rrVert_{P}^2
\bigr)
\\
&&\quad\leq(1+\eta)\llVert P_{i|S}-P_{i|V}\rrVert
_{P}^2+r_2\llVert\Ph_{i|V}-P_{i|V}
\rrVert_{P}^2
\\
&&\quad\quad{}+(1+r_2) \bigl( \llVert\Ph_{i|V}-P_{i|V}\rrVert
_{\Ph}^2-\llVert\Ph_{i|V}-P_{i|V}
\rrVert_{P}^2 \bigr)+ \biggl(\frac{17}{\eta
}+1 \biggr)
\frac
{\log
(N_s^2\delta)}{3n}.
\end{eqnarray*}
Let $C$ be the constant given by Lemma B.5 and let
\[
\Omega_{*}= \bigl\{\forall V\in\Vc_{s}, \bigl| \llVert \Ph
_{i|V}-P_{i|V}\rrVert_{\Ph}^2-\llVert \Ph_{i|V}-P_{i|V}\rrVert
_{P}^2 \bigr|\leq C
\varepsilon\llVert\Ph_{i|V}-P_{i|V}\rrVert_{P}^2
\bigr\}.
\]
It comes from Lemma B.5 that $P(\Omega_{*})\geq1-\delta
^{-1}$. Moreover, on $\Omega_{\mathrm{comp}}\cap\Omega_{*}$, we have, from
(\ref{eq:MinimCrit}), for all $0<\eta<1$,
\begin{eqnarray*}
&&(1-\eta)\llVert P_{i|S}-P_{i|\Vh}\rrVert_{P}^2+
\bigl(r_1-C(1+r_1)\varepsilon\bigr)\llVert
\Ph_{i|\Vh}-P_{i|\Vh}\rrVert_{P}^2
\\
&&\quad\leq(1+\eta)\llVert P_{i|S}-P_{i|V}\rrVert
_{P}^2+\bigl(r_2+C(1+r_2)
\varepsilon\bigr)\llVert\Ph_{i|V}-P_{i|V}\rrVert
_{P}^2+\frac{6}{\eta}\frac{\log(N_s^2\delta)}{n}.
\end{eqnarray*}

\end{appendix}

\section*{Acknowledgements}
We are grateful to Antonio Galves for many discussions
and fruitful advices during the redaction of the paper. We
also would like to thank the referees for the comments that
considerably
improved the manuscript.

ML was supported by FAPESP Grant 2009/09494-0. DYT was partially
supported by FAPESP Grant 2008/08171-0, Pew Latin American Fellowship,
and Ci\^encia sem Fronteiras Fellowship (CNPq Grant  246778/2012-1).
This work is part of USP project ``Mathematics, computation, language
and the brain''.

\begin{supplement}%[id=suppA]
\stitle{Supplement to ``Sharp oracle inequalities and slope heuristic
for specification~probabilities estimation
in discrete random fields''\\}
\slink[doi]{10.3150/14-BEJ660SUPP} %[doi,text={...}] - jei reikia
%suskaldyti doi
\sdatatype{.pdf}
\sfilename{BEJ660\_supp.pdf}
\sdescription{On this supplementary material available on-line, we
prove the probabilistic tools needed in the proofs of the main results.
The second part provides additional simulation results. The last one is
devoted to the extension of all our results to the K\"ullback loss.}
\end{supplement}

% imsref loaded by arune.pranskunaite, 2014-07-25 14:47:47
%
% imsref loaded by arune.pranskunaite, 2014-07-28 10:02:19

\printhistory

\begin{thebibliography}{25}

%b1 #&#
\bibitem{AB10}
%
\begin{bincollection}[auto:STB|2014/06/18|12:29:53]
\bauthor{\bsnm{Arlot},~\bfnm{S.}\binits{S.}} \AND
\bauthor{\bsnm{Bach},~\bfnm{F.}\binits{F.}}
(\byear{2010}).
\btitle{Data-driven calibration of linear estimators with minimal penalties}.
In \bbooktitle{Advances in Neural Information Processing Systems (NIPS)}
(\beditor{\bfnm{Y.}\binits{Y.}~\bsnm{Bengio}},
\beditor{\bfnm{D.}\binits{D.}~\bsnm{Schuurmans}},
\beditor{\bfnm{J.~D.}\binits{J.D.}~\bsnm{Lafferty}},
\beditor{\bfnm{C.~K.~I.}\binits{C.K.I.}~\bsnm{Williams}} \AND
\beditor{\bfnm{A.}\binits{A.}~\bsnm{Culotta}}, eds.)
\bvolume{22}
\bpages{46--54}.
\bnote{Available at \surl{http://papers.nips.cc/book/\\advances-in-neural-information-processing-systems-22-2009}}.
\end{bincollection}
%
\bptok{imsref}%
\endbibitem


%b2 #&#
\bibitem{AM08}
%
\begin{barticle}[auto:STB|2014/06/18|12:29:53]
\bauthor{\bsnm{Arlot},~\bfnm{S.}\binits{S.}} \AND
\bauthor{\bsnm{Massart},~\bfnm{P.}\binits{P.}}
(\byear{2009}).
\btitle{Data-driven calibration of penalties for least-squares regression}.
\bjournal{J. Mach. Learn. Res.}
\bvolume{10}
\bpages{245--279}.
\end{barticle}
%
\bptok{imsref}%
\endbibitem

%b3 #&#
\bibitem{BBM99}
%
\begin{barticle}[mr]
\bauthor{\bsnm{Barron},~\bfnm{Andrew}\binits{A.}},
\bauthor{\bsnm{Birg{\'e}},~\bfnm{Lucien}\binits{L.}} \AND
\bauthor{\bsnm{Massart},~\bfnm{Pascal}\binits{P.}}
(\byear{1999}).
\btitle{Risk bounds for model selection via penalization}.
\bjournal{Probab. Theory Related Fields}
\bvolume{113}
\bpages{301--413}.
\bid{doi={10.1007/s004400050210}, issn={0178-8051}, mr={1679028}}
\end{barticle}
%
\bptok{imsref}%
% NOT OUTPUTED:
% issn = 0178-8051
% url = http://dx.doi.org/10.1007/s004400050210
% number = 3
% coden = PTRFEU
% fjournal = Probability Theory and Related Fields
\endbibitem

%b4 #&#
\bibitem{BS91}
%
\begin{barticle}[mr]
\bauthor{\bsnm{Barron},~\bfnm{Andrew~R.}\binits{A.R.}} \AND
\bauthor{\bsnm{Sheu},~\bfnm{Chyong-Hwa}\binits{C.-H.}}
(\byear{1991}).
\btitle{Approximation of density functions by sequences of exponential
families}.
\bjournal{Ann. Statist.}
\bvolume{19}
\bpages{1347--1369}.
\bid{doi={10.1214/aos/1176348252}, issn={0090-5364}, mr={1126328}}
\end{barticle}
%
\bptok{imsref}%
% NOT OUTPUTED:
% issn = 0090-5364
% url = http://dx.doi.org/10.1214/aos/1176348252
% number = 3
% coden = ASTSC7
% fjournal = The Annals of Statistics
\endbibitem

%b5 #&#
\bibitem{Bento09}
%
\begin{bmisc}[auto:STB|2014/06/18|12:29:53]
\bauthor{\bsnm{Bento},~\bfnm{J.}\binits{J.}} \AND
\bauthor{\bsnm{Montanari},~\bfnm{A.}\binits{A.}} (\byear{2009}).
\bhowpublished{Which graphical models are difficult to learn?
Available at \surl{http://arxiv.org/pdf/0910.5761}}.
\end{bmisc}
%
\bptok{imsref}%
% NOT OUTPUTED:
% sortkey = Bento(2009
\endbibitem

%b6 #&#
\bibitem{BM97}
%
\begin{bincollection}[mr]
\bauthor{\bsnm{Birg{\'e}},~\bfnm{Lucien}\binits{L.}} \AND
\bauthor{\bsnm{Massart},~\bfnm{Pascal}\binits{P.}}
(\byear{1997}).
\btitle{From model selection to adaptive estimation}.
In \bbooktitle{Festschrift for {L}ucien {L}e {C}am}
\bpages{55--87}.
\blocation{New York}:
\bpublisher{Springer}.
\bid{mr={1462939}}
\end{bincollection}
%
\bptok{imsref}%
\endbibitem

%b7 #&#
\bibitem{BM01}
%
\begin{barticle}[mr]
\bauthor{\bsnm{Birg{\'e}},~\bfnm{Lucien}\binits{L.}} \AND
\bauthor{\bsnm{Massart},~\bfnm{Pascal}\binits{P.}}
(\byear{2001}).
\btitle{Gaussian model selection}.
\bjournal{J. Eur. Math. Soc. (JEMS)}
\bvolume{3}
\bpages{203--268}.
\bid{doi={10.1007/s100970100031}, issn={1435-9855}, mr={1848946}}
\end{barticle}
%
\bptok{imsref}%
% NOT OUTPUTED:
% issn = 1435-9855
% url = http://dx.doi.org/10.1007/s100970100031
% number = 3
% fjournal = Journal of the European Mathematical Society (JEMS)
\endbibitem

%b8 #&#
\bibitem{BM07}
%
\begin{barticle}[mr]
\bauthor{\bsnm{Birg{\'e}},~\bfnm{Lucien}\binits{L.}} \AND
\bauthor{\bsnm{Massart},~\bfnm{Pascal}\binits{P.}}
(\byear{2007}).
\btitle{Minimal penalties for {G}aussian model selection}.
\bjournal{Probab. Theory Related Fields}
\bvolume{138}
\bpages{33--73}.
\bid{doi={10.1007/s00440-006-0011-8}, issn={0178-8051}, mr={2288064}}
\end{barticle}
%
\bptok{imsref}%
% NOT OUTPUTED:
% issn = 0178-8051
% url = http://dx.doi.org/10.1007/s00440-006-0011-8
% number = 1-2
% coden = PTRFEU
% fjournal = Probability Theory and Related Fields
\endbibitem

%b9 #&#
\bibitem{Bo02}
%
\begin{barticle}[mr]
\bauthor{\bsnm{Bousquet},~\bfnm{Olivier}\binits{O.}}
(\byear{2002}).
\btitle{A {B}ennett concentration inequality and its application to
suprema of empirical processes}.
\bjournal{C. R. Math. Acad. Sci. Paris}
\bvolume{334}
\bpages{495--500}.
\bid{doi={10.1016/S1631-073X(02)02292-6}, issn={1631-073X}, mr={1890640}}
\end{barticle}
%
\bptok{imsref}%
% NOT OUTPUTED:
% issn = 1631-073X
% url = http://dx.doi.org/10.1016/S1631-073X(02)02292-6
% number = 6
% fjournal = Comptes Rendus Math\'ematique. Acad\'emie des Sciences.
%Paris
\endbibitem

%b10 #&#
\bibitem{Bresler08}
%
\begin{bincollection}[mr]
\bauthor{\bsnm{Bresler},~\bfnm{Guy}\binits{G.}},
\bauthor{\bsnm{Mossel},~\bfnm{Elchanan}\binits{E.}} \AND
\bauthor{\bsnm{Sly},~\bfnm{Allan}\binits{A.}}
(\byear{2008}).
\btitle{Reconstruction of {M}arkov random fields from samples: Some
observations and algorithms}.
In \bbooktitle{Approximation, Randomization and Combinatorial Optimization}.
\bseries{Lecture Notes in Computer Science}
\bvolume{5171}
\bpages{343--356}.
\blocation{Berlin}:
\bpublisher{Springer}.
\bid{doi={10.1007/978-3-540-85363-3_28}, mr={2538799}}
\end{bincollection}
%
\bptok{imsref}%
% NOT OUTPUTED:
% url = http://dx.doi.org/10.1007/978-3-540-85363-3_28
\endbibitem

%b11 #&#
\bibitem{Brown04}
%
\begin{barticle}[auto:STB|2014/06/18|12:29:53]
\bauthor{\bsnm{Brown},~\bfnm{E.~N.}\binits{E.N.}},
\bauthor{\bsnm{Kass},~\bfnm{R.~E.}\binits{R.E.}} \AND
\bauthor{\bsnm{Mitra},~\bfnm{P.~P.}\binits{P.P.}}
(\byear{2004}).
\btitle{Multiple neural spike train data analysis: State-of-the-art and
future challenges}.
\bjournal{Nature Neuroscience}
\bvolume{7}
\bpages{456--461}.
\end{barticle}
%
\bptok{imsref}%
\endbibitem

%b12 #&#
\bibitem{Csiszar06}
%
\begin{barticle}[mr]
\bauthor{\bsnm{Csisz{\'a}r},~\bfnm{Imre}\binits{I.}} \AND
\bauthor{\bsnm{Talata},~\bfnm{Zsolt}\binits{Z.}}
(\byear{2006}).
\btitle{Consistent estimation of the basic neighborhood of {M}arkov
random fields}.
\bjournal{Ann. Statist.}
\bvolume{34}
\bpages{123--145}.
\bid{doi={10.1214/009053605000000912}, issn={0090-5364}, mr={2275237}}
\end{barticle}
%
\bptok{imsref}%
% NOT OUTPUTED:
% issn = 0090-5364
% url = http://dx.doi.org/10.1214/009053605000000912
% number = 1
% coden = ASTSC7
% fjournal = The Annals of Statistics
\endbibitem

%b13 #&#
\bibitem{CT06}
%
\begin{barticle}[mr]
\bauthor{\bsnm{Csisz{\'a}r},~\bfnm{Imre}\binits{I.}} \AND
\bauthor{\bsnm{Talata},~\bfnm{Zsolt}\binits{Z.}}
(\byear{2006}).
\btitle{Context tree estimation for not necessarily finite memory
processes, via BIC and {MDL}}.
\bjournal{IEEE Trans. Inform. Theory}
\bvolume{52}
\bpages{1007--1016}.
\bid{doi={10.1109/TIT.2005.864431}, issn={0018-9448}, mr={2238067}}
\end{barticle}
%
\bptok{imsref}%
% NOT OUTPUTED:
% issn = 0018-9448
% url = http://dx.doi.org/10.1109/TIT.2005.864431
% number = 3
% coden = IETTAW
% fjournal = Institute of Electrical and Electronics Engineers.
%Transactions on Information Theory
\endbibitem

%b14 #&#
\bibitem{GOT10}
%
\begin{bmisc}[auto:STB|2014/06/18|12:29:53]
\bauthor{\bsnm{Galves},~\bfnm{A.}\binits{A.}},
\bauthor{\bsnm{Orlandi},~\bfnm{E.}\binits{E.}} \AND
\bauthor{\bsnm{Takahashi},~\bfnm{D.~Y.}\binits{D.Y.}} (\byear{2010}).
\bhowpublished{Identifying interacting pairs of sites in
infinite range ising models. Preprint. Available at
\surl{http://arxiv.org/\\abs/1006.0272}}.
\end{bmisc}
%
\bptok{imsref}%
% NOT OUTPUTED:
% sortkey = Galves(2010
\endbibitem

%b15 #&#
\bibitem{Giorgii88}
%
\begin{bbook}[mr]
\bauthor{\bsnm{Georgii},~\bfnm{Hans-Otto}\binits{H.-O.}}
(\byear{1988}).
\btitle{Gibbs Measures and Phase Transitions}.
\bseries{de Gruyter Studies in Mathematics}
\bvolume{9}.
\blocation{Berlin}:
\bpublisher{de Gruyter}.
\bid{doi={10.1515/9783110850147}, mr={0956646}}
\end{bbook}
%
\bptok{imsref}%
% NOT OUTPUTED:
% isbn = 0-89925-462-4
% url = http://dx.doi.org/10.1515/9783110850147
% fpage = xiv+525
\endbibitem

%b16 #&#
\bibitem{Ler2010mixing}
%
\begin{barticle}[mr]
\bauthor{\bsnm{Lerasle},~\bfnm{Matthieu}\binits{M.}}
(\byear{2011}).
\btitle{Optimal model selection for density estimation of stationary
data under various mixing conditions}.
\bjournal{Ann. Statist.}
\bvolume{39}
\bpages{1852--1877}.
\bid{doi={10.1214/11-AOS888}, issn={0090-5364}, mr={2893855}}
\end{barticle}
%
\bptok{imsref}%
% NOT OUTPUTED:
% issn = 0090-5364
% url = http://dx.doi.org/10.1214/11-AOS888
% number = 4
% fjournal = The Annals of Statistics
\endbibitem

%b17 #&#
\bibitem{Le09}
%
\begin{barticle}[mr]
\bauthor{\bsnm{Lerasle},~\bfnm{Matthieu}\binits{M.}}
(\byear{2012}).
\btitle{Optimal model selection in density estimation}.
\bjournal{Ann. Inst. Henri Poincar\'e Probab. Stat.}
\bvolume{48}
\bpages{884--908}.
\bid{doi={10.1214/11-AIHP425}, issn={0246-0203}, mr={2976568}}
\end{barticle}
%
\bptok{imsref}%
% NOT OUTPUTED:
% issn = 0246-0203
% url = http://dx.doi.org/10.1214/11-AIHP425
% number = 3
% fjournal = Annales de l'Institut Henri Poincar\'e Probabilit\'es et
%Statistiques
\endbibitem

%b18 #&#
\bibitem{LT10}
%
\begin{barticle}[mr]
\bauthor{\bsnm{Lerasle},~\bfnm{Matthieu}\binits{M.}} \AND
\bauthor{\bsnm{Takahashi},~\bfnm{Daniel~Y.}\binits{D.Y.}}
(\byear{2011}).
\btitle{An oracle approach for interaction neighborhood estimation
in random fields}.
\bjournal{Electron. J. Stat.}
\bvolume{5}
\bpages{534--571}.
\bid{doi={10.1214/11-EJS618}, issn={1935-7524}, mr={2813554}}
\end{barticle}
%
\bptok{imsref}%
% NOT OUTPUTED:
% issn = 1935-7524
% url = http://dx.doi.org/10.1214/11-EJS618
% fjournal = Electronic Journal of Statistics
\endbibitem

%b25 #&#
\bibitem{supp}
\begin{bmisc}[author]
\bauthor{\bsnm{Lerasle}~~\bfnm{Matthieu}\binits{M.}} \and
\bauthor{\bsnm{Takahashi}~~\bfnm{Daniel~Y.}\binits{D. Y.}}
(\byear{2014}).
\bhowpublished{Supplement to ``Sharp oracle inequalities
and slope heuristic for specification probabilities estimation
in discrete random fields.''
DOI:\doiurl{10.3150/14-BEJ660SUPP}}.
\bptok{imsref}\end{bmisc}
\endbibitem

%b19 #&#
\bibitem{Ma07}
%
\begin{bbook}[mr]
\bauthor{\bsnm{Massart},~\bfnm{Pascal}\binits{P.}}
(\byear{2007}).
\btitle{Concentration Inequalities and Model Selection}.
\bseries{Lecture Notes in Math.}
\bvolume{1896}.
\blocation{Berlin}:
\bpublisher{Springer}.
\bnote{Lectures from the 33rd Summer School on Probability Theory held
in Saint-Flour, July 6--23, 2003. With a foreword by Jean Picard}.
\bid{mr={2319879}}
\end{bbook}
%
\bptok{imsref}%
% NOT OUTPUTED:
% isbn = 978-3-540-48497-4; 3-540-48497-3
% fpage = xiv+337
\endbibitem

%b20 #&#
\bibitem{Mizuseki09}
%
\begin{bmisc}[auto:STB|2014/06/18|12:29:53]
\bauthor{\bsnm{Pastalkova},~\bfnm{E.}\binits{E.}},
\bauthor{\bsnm{Buzs\'aki},~\bfnm{G.}\binits{G.}},
\bauthor{\bsnm{Mizuseki},~\bfnm{K.}\binits{K.}} \AND
\bauthor{\bsnm{Sirota},~\bfnm{A.}\binits{A.}}
\bhowpublished{Theta oscillations provide temporal windows for local
circuit
computation in the entorhinal-hippocampal loop.
{\it Neuron} \textbf{64} 267--280}.
\end{bmisc}
%
\bptok{imsref}%
% NOT OUTPUTED:
% sortkey = E
\endbibitem

%b21 #&#
\bibitem{PW10}
%
\begin{barticle}[mr]
\bauthor{\bsnm{Ravikumar},~\bfnm{Pradeep}\binits{P.}},
\bauthor{\bsnm{Wainwright},~\bfnm{Martin~J.}\binits{M.J.}} \AND
\bauthor{\bsnm{Lafferty},~\bfnm{John~D.}\binits{J.D.}}
(\byear{2010}).
\btitle{High-dimensional {I}sing model selection using {$\ell\sb
1$}-regularized logistic regression}.
\bjournal{Ann. Statist.}
\bvolume{38}
\bpages{1287--1319}.
\bid{doi={10.1214/09-AOS691}, issn={0090-5364}, mr={2662343}}
\end{barticle}
%
\bptok{imsref}%
% NOT OUTPUTED:
% issn = 0090-5364
% url = http://dx.doi.org/10.1214/09-AOS691
% number = 3
% coden = ASTSC7
% fjournal = The Annals of Statistics
\endbibitem

%b22 #&#
\bibitem{Sa13}
%
\begin{barticle}[auto:STB|2014/06/18|12:29:53]
\bauthor{\bsnm{Saumard},~\bfnm{A.}\binits{A.}}
(\byear{2013}).
\btitle{The slope heuristics in heteroscedastic regression}.
\bjournal{Electron. J. Stat.}
\bvolume{7}
\bpages{1184--1223}.
\bid{mr={3056072}}
\end{barticle}
%
\bptok{imsref}%
% NOT OUTPUTED:
% sortkey = Saumard(2013
\endbibitem

%b23 #&#
\bibitem{Schneidman06}
%
\begin{barticle}[auto:STB|2014/06/18|12:29:53]
\bauthor{\bsnm{Schneidman},~\bfnm{E.}\binits{E.}},
\bauthor{\bsnm{Berry},~\bfnm{M.~J.}\binits{M.J.}},
\bauthor{\bsnm{Segev},~\bfnm{R.}\binits{R.}} \AND
\bauthor{\bsnm{Bialek},~\bfnm{W.}\binits{W.}}
(\byear{2006}).
\btitle{Weak pairwise correlations imply strongly correlated network
states in a neural population}.
\bjournal{Nature}
\bvolume{440}
\bpages{1007--1012}.
\end{barticle}
%
\bptok{imsref}%
\endbibitem

%b24 #&#
\bibitem{Takahashi10}
%
\begin{barticle}[auto:STB|2014/06/18|12:29:53]
\bauthor{\bsnm{Takahashi},~\bfnm{N.}\binits{N.}},
\bauthor{\bsnm{Sasaki},~\bfnm{T.}\binits{T.}},
\bauthor{\bsnm{Matsumoto},~\bfnm{W.}\binits{W.}} \AND
\bauthor{\bsnm{Ikegaya},~\bfnm{Y.}\binits{Y.}}
(\byear{2010}).
\btitle{Circuit topology for synchronizing neurons in spontaneously
active networks}.
\bjournal{Proc. Natl. Acad. Sci. USA}
\bvolume{107}
\bpages{10244--10249}.
\end{barticle}
%
\bptok{imsref}%
\endbibitem



\end{thebibliography}
\end{document}